\documentclass[13pt. oneside]{amsart} 
\usepackage[lmargin=1in,rmargin=1in,
bmargin=1.1in, tmargin=1.1in]{geometry}
\usepackage[dvipsnames]{xcolor}
\usepackage{amsmath, amssymb,tikz}
\usepackage{tikz-cd}
\usepackage{tikz}
\usepackage{todonotes}
\usepackage{mathrsfs}
\usepackage{extarrows}
\usepackage{graphicx}
\usepackage[breaklinks, pagebackref]{hyperref}
\usepackage{IEEEtrantools}
\usepackage{mathrsfs}
\usepackage[utf8]{inputenc}
\usepackage[english]{babel}
\usepackage{comment}
\usepackage{csquotes}
\usepackage[shortlabels]{enumitem}
\usepackage{stmaryrd}
\usepackage{manfnt}

\DeclareSymbolFont{bbold}{U}{bbold}{m}{n}
\DeclareSymbolFontAlphabet{\mathbbold}{bbold}

\hypersetup{
    colorlinks=true,
    linkcolor=blue,
    filecolor=magenta, 
    citecolor=red,     
    urlcolor=cyan,
    }

\def\im{\mathrm{im}}
\def\diag{\mathrm{diag}}

\def\Ch{\mathrm{Ch}}

\def\Ann{\mathrm{Ann}}

\def\Ann{\mathrm{Ann}}                                                               
\usepackage{amsmath,amsthm,amssymb}
\usepackage{colonequals}

\newcommand{\ZZ}{\mathbb{Z}}

\newcounter{dummypart}



\numberwithin{equation}{section}

\newtheorem{theorem}[equation]{Theorem}
\newtheorem{proposition}[equation]{Proposition}
\newtheorem{lemma}[equation]{Lemma}
\newtheorem{corollary}[equation]{Corollary}
\newtheorem*{corollary*}{Corollary}
\newtheorem{conjecture}[equation]{Conjecture}

\theoremstyle{definition}

\newtheorem*{notation*}{Notation}

\theoremstyle{definition}

\newtheorem{definition}[equation]{Definition}

\theoremstyle{remark}
\newtheorem{remark}[equation]{Remark}

\newtheorem{note*}[equation]{Note}

\tikzcdset{scale cd/.style={every label/.append style={scale=#1},
    cells={nodes={scale=#1}}}}

\usepackage{tikz-cd}
\usepackage{mathrsfs}

\DeclareFontFamily{U}{wncy}{}
\DeclareFontShape{U}{wncy}{m}{n}{<->wncyr10}{}
\DeclareSymbolFont{mcy}{U}{wncy}{m}{n}
\DeclareMathSymbol{\sha}{\mathord}{mcy}{"58}

\makeatletter
\renewcommand*\env@matrix[1][\arraystretch]{%
  \edef\arraystretch{#1}%
  \hskip -\arraycolsep
  \let\@ifnextchar\new@ifnextchar
  \array{*\c@MaxMatrixCols c}}
\makeatother

\usepackage{mathtools,bm}

\DeclarePairedDelimiterX\Set[1]\{\}{%
  #1%
}

\usepackage{stmaryrd}

\newcommand{\GL}{\mathrm{GL}}

\begin{document}

\title[Smith Normal Form of Matrices Associated with differential posets]{Smith normal form of matrices associated with differential Posets}

\author{Syed Waqar Ali Shah}

\address{Department of Mathematics\\
Lahore University of Management Sciences\\
Lahore, Pakistan}

\email{swaqar.66080@gmail.com}

\begin{abstract}
We prove a conjecture of Miller and Reiner on the existence of Smith normal form for the  $DU$-operators for a certain class of $r$-differential posets.
\end{abstract}



\maketitle

\tableofcontents

\section{Introduction}           

    Let $r$ be a positive integer. We say that a partially ordered set (poset) $ P $ is \emph{$ r$-differential}, if it satisfies the following three conditions:    
\begin{enumerate} 
\item[(D1)]$P$ is graded, locally finite, has all ranks finite and possesses a unique minimal element.    
\item[(D2)]If two distinct elements of $P$ have exactly $k$ elements
that are covered by both of them, then there are  exactly $k$ elements
that cover them both.      
\item[(D3)] If an element of $P$ covers exactly  $k$ elements, then it is covered by exactly $k+r$  elements.                                             
\end{enumerate}
Associated to every $ r $-differential poset are two families of maps, known as \textit{up} and \textit{down} maps. Let $ P_{n} $ be the $ n$-th rank of $ P $, which we take to be the empty set if $ n < 0 $, and set $ p_{n} : = | P_{n} |$. For any commutative ring $ R$ with identity and characteristic $ 0 $, let $ 
RP_{n} \cong R^{p_{n}}$ be the free module over $ R $ with basis $ P_{n} $. We define 
$$
\begin{aligned}
U_n : RP_{n} & \rightarrow RP_{n+1}\\
D_n :  RP_{n}  & \rightarrow RP_{n-1}
\end{aligned}
$$
for all $ n \geq 0 $ on basis elements as follows:  $ U_{n} $ sends  $ x \in P_{n} $ to the sum (with coefficients $1$) of all elements in $ P_{n+1} $ that cover $ x $ and $ D_{n} $ sends $ x $ to the sum of all elements in $ P_{n-1} $ that are covered by $ x $.  We then define
$$
\begin{aligned}
UD_n &:=U_{n-1} \circ D_n \\
DU_n &:=D_{n+1} \circ U_n
\end{aligned}
$$
The two conditions (D2) and (D3) can then be recast as 
$$    DU_{n}  -   UD_{n} = r \cdot \mathbbold{1} . $$ 
The most well-known examples of $ 1 $-differential posets are the Young's lattice $ \mathbf{Y} $ and Young-Fibonacci lattice $  \mathbf{Y}F $.  Their $ r $-fold cartesian products are examples of  $r$-differential posets. Another important example of an $ r$-differential poset is the $ r $-Fibonacci poset $ Z(r) $ where   $ Z(1) = \mathbf{Y} F$.

Differential posets were first defined by Stanley in \cite{Stanley88} with up and down maps defined over fields. Later, Miller and Reiner defined them over arbitrary rings in \cite{Miller} and conjectured a remarkable property of the $ DU_{n}    $-operators over the ring of integers which we now describe.  Recall that a $ m \times n$ matrix $ A = (a_{ij}) $ over an integral domain $ R $ is said to have \textit{Smith normal form} (SNF) over $ R $ if there exist invertible matrices $ P \in R ^{ m \times m}$, $ Q \in R ^{ n \times n} $ such that $ B = PAQ $ is a diagonal matrix in the sense that $ b_{ij} = 0 $ if $ i \neq j $, and $ s_{i} := b_{ii} $ for $ 1 \leq i \leq k =   \min  \left  \{     m ,  n     \right \} $ satisfy the divisibilities $$ s_{1}| s_{2} | \ldots | s_{k} . $$  It is known that  if     $ R $ is a PID, any matrix $ A $ always has  an   SNF which is unique in the sense that the diagonal entries $ s_{i} $ are unique up to units of $ R $. If $ R $ is not a PID, such a form does not necessarily exist. If it does exist however, it is  still   unique.   

Assume now that $ R =  \ZZ $. Let $ [DU_{n}]$ be the matrix of $ DU_{n} $ with respect to the standard basis of $   \ZZ P_{n} $ and $ I_{p_{n}} $ be the $ p_{n} \times p_{n} $ identity matrix.    
\begin{conjecture}{\cite[Miller-Reiner]{Miller}}  For all $n $, the matrix $  [DU_{n}] + xI_{p_{n}} $ has SNF over $ \ZZ[x] $.      
\label{mainconj}                                                                       
\end{conjecture}
Miller and Reiner verified this conjecture for the $ r $-Fibonacci poset $Z(r)$ in \cite{Miller}. The problem was later investigated by Cai and Stanley in \cite{SNFforYoung} for the case $ \mathbf{Y}^{r} $ and the case $ r = 1 $ was settled in the affirmative. As noted in the survey \cite{SNFsurvey}, the case $ r > 1 $ was later handled by Zipei Nie, though  a written proof has not yet  appeared.                               

In this paper, we prove this conjecture for any $ r$-differential poset that satisfies certain conditions which are stated in Theorem \ref{mainteo}. These conditions 
are  closely related to  two additional conjectures (Conjectures $ 2.3 $ and $ 2.4$ of \cite{Miller}) made by Miller and Reiner.    Our strategy is to study the $ \mathbb{Z}[x] $-module structure of $ \mathbb{Z} P_{n}  $ where the action of $ x $   is induced by the operator $ DU_{n}    $. The statement on the Smith normal form of $ [DU_{n}] + x I_{p_{n}} $ is then easily translated into a statement on the existence of \emph{integral canonical form} for $DU_{n} $ which   is   an integral analogue of the rational canonical form.  This allows us to relate the structure of $ DU_{n} $ and $ DU_{n+1} $    which     paves the way for an induction argument. 
\bigskip

\noindent    \textbf{Acknowledgements:} I would like to specially thank my mentor, Dr. Shaheen Nazir, for introducing me to this topic and providing help and guidance throughout my research journey and Ivan Sadofschi Costa, for providing help with various computations in Sage, helpful edits and for being a constant source of support. I would also like to thank Brian Conrad for his valuable feedback and insightful conversations. In addition, I thank Victor Reiner for his comments on an earlier version, Alexander Miller for sharing the data on rank sizes and Richard Stanley for useful  communication. Finally, I wish to thank the anonymous referees for pointing out an error and for their valuable comments and suggestions.
\section{Recollections}   
\label{recoll}

We start with a theorem of Stanley. 
\begin{theorem} \cite[\S 4]{Stanley88}    \label{Stanley}         Let $P$ be  an  $r$-differential poset and $ R $ a field of characteristic $ 0 $. Then
\[
\begin{aligned}
\Ch(DU_{n}) &=\prod_{j=0}^n (x - r j - r  )^{\Delta p_{n-j}} \\
\Ch(UD_n) &=\prod_{j=0}^n (x - r j  )^{\Delta p_{n-j}}
\end{aligned}
\]
where $ \Ch(A) = \Ch(A,x) $ denotes the characteristic polynomial of the operator $ A$, and $ \Delta p_n:=p_n - p_{n-1} $ denotes the rank difference.  Furthermore, the operators $ DU_{n} $ and $ UD_{n} $ are diagonalizable.                                                                                   
\end{theorem}

We make some immediate conclusions. First, the rank function is non-decreasing as $ \Delta   p_{n} $ must be non-negative in the expressions above. Second, $ DU_{n} $ is invertible as all its eigenvalues are non-zero.  Thus $ U_{n} $ is injective and $ D_{n+1} $ is surjective for all $ n \geq 0 $.   Third, the  sequence of  invariant factors of the torsion  $ R[x] $-module $ RP_{n} $ where $ x $ acts via $ DU_{n} $ is uniquely determined by $  \Ch (DU_{n}) $ for each  $ n  $.  More precisely, there exists a decomposition  $$ RP_{n} =  M_{1} \oplus M_{2} \oplus  \ldots \oplus  M_{k} $$ where $  k  : =  \max           \left \{    \Delta p_{0}, \ldots, \Delta p_{n} \right \}   $  and each $ M_{i} $ is a monogenic  torsion $R[x] $-module with annihilator  \begin{equation} 
\label{exp}  a_{i}(x)  =   \prod _{       \substack{ j \in \left  \{  0 , \ldots, n \right \}    \\ \Delta p_{n-j}  \geq k - i   +  1 }  }  ( x - rj- r)      .     
\end{equation}  Note that $  a_{1} (x) | a_{2}(x) | \ldots  |  a_{m  } ( x)    $. That $ a_{i}(x) $ must be given by (\ref{exp})     follows from diagonalizability of $ DU_{n} $ which forces $ a_{i}(x) $ to have no repeated  factors.   Since the non-constant polynomials in the Smith normal form of  $ xI_{p_{n}} -   [DU_{n}]  
$  are the invariant factors of $ R P_{n}   $ as an $R[x]$-module, the SNF of $ xI_{p_{n}} - [DU_{n}]    $     over $ R[x] $ is  $$ \diag( 1, \ldots, 1, a_{1}(x) , a_{2}(x), \ldots, a_{m}(x) ) . $$          
See \cite[\S 8.2]{Miller} for  details.

Many of  these conclusions fail when one replaces the field 
$ R $ by  a PID of characteristic zero. For instance, $  DU_{n} $ and $ UD_{n} $ are not necessarily diagonalizable over $ R $.  An explicit  counterexample is given by $ DU_{2} $ for $ P         =    \mathbf{Y}   $. In the standard basis the matrix  is $ A = \left ( \begin{smallmatrix} 2  & 1  \\ 1  &  2          \end{smallmatrix}   \right ) $.   Diagonalizability of $ A $ is equivalent to requiring  that $ R^{2} =  \ker (DU_{2} - I_{2}) \oplus \ker ( DU_{2} - 3I_{2}) $ or that  
$$    R  =   R \left ( \begin{smallmatrix}    1 \\  - 1 \end{smallmatrix} \right ) \oplus  R \left ( \begin{smallmatrix} 1 \\  1  \end{smallmatrix} \right  )  $$
This is however absurd if $ 2 \notin R^{\times }$.  Invertibility of $ DU_{n} $ is also no longer guaranteed as the determinant of its matrix may not land in $ R^{\times} $.  For instance,  $ \det(A) = 3 $ in the example above. In particular, $ DU_{n} $ is not necessarily surjective and one cannot deduce the surjectivity of down maps.   Finally since $ R[x] $ is not a PID, one in general does not expect a similar decomposition to exist for $  R P_{n} $ 
as an $ R[x]$-module.  However, the existence of such a decomposition is equivalent to the conjecture of Miller and Reiner. We explain this in the next section.

\section{Integral canonical forms} 

Throughout this section,  $ R $ denotes an integral domain and $ F $ its field of fractions. For an $ R[x] $-module $ M $, we let  $ \mathrm{Ann}(M) $ denote  the annihilator of $ M $. For a monic polynomial $ a(x) \in R[x] $, we denote its Frobenius companion matrix by $  \mathcal{C}_{a(x)}  $. For $ m , n $ positive integers, we let $ \mathrm{Mat}_{m \times n }(R) $   denote the set of all $ m \times n $-matrices with entires in $ R $. Two matrices $ A , B \in \mathrm{Mat}_{m \times n }(R[x]) $ are said to be congruent if there exist $ P \in \GL_{m}(R[x]) $, $ Q \in \GL_{n}(R[x]) $ such that $ B =  P A Q $. 

\begin{definition} An $ R[x]$-module is  said to have an \textit{invariant factor decomposition} if there exist monogenic $ R[x]$-submodules   $ M_{1} ,  M_{2} ,  \ldots,  M_{k}   $ of $ M  $  such that  $ \mathrm{Ann}(M_{i}) $ is generated by a monic, non-constant polynomial $  a_{i} = a_{i}(x) $ that satisfy the chain of divisibilities $$  a_{1} |  a_{2}  |            \cdots  |  a_{k}       $$            
and  $  M    =   M_{1}     \oplus   M_{2}   \oplus    \cdots   \oplus      M_{k}            $. 
\end{definition} 

If  $ M $ is an $R[x]$-module admitting such a decomposition, it is necessarily free of finite rank over $ R $.   Tensoring  both  sides       of     $ M = M_{1} \oplus  \cdots. \oplus   M_{k} $ with $ F $ yields the usual invariant factor decomposition of the $ F[x] $-module $ M \otimes_{R} F $. The uniqueness of invariant factor decomposition over fields implies that the sequence  of   polynomials $ a_{i}  $ is uniquely determined by $ M  $. We will therefore  refer to $ a_{1}, \ldots, a_{k} $ as  the \emph{invariant factor sequence}  for $ M $.  
 
Let us denote  $    d_{i}  : =     \deg   a_{i}(x)  $, $ v_{i} \in M $ an $ R[x]$-generator and $ \varphi : M \to M $ the $ R $-linear endomorphism given by multiplication by $ x $. Then $$           \alpha            : =   ( v_{1} ,  \ldots,  x ^ { d_{1}    - 1 } v_{1} ,  v_{2}  , \ldots,   x ^{ d_{2} - 1}  v_{2} ,  \ldots ,   v_{k}  , \ldots ,         x ^ { d_{k} -1 } v_{k}  )   $$ is an $ R $-basis for $  M  $ with respect to which the matrix of $ \varphi $ is block-diagonal with $ k $ block matrices, the $  i $-th block being the companion matrix of $ a_{i}(x) $. In other words,                                                                   
\begin{equation} \label{xalpha}    [\varphi]_{\alpha}  = \begin{pmatrix}    \mathcal{C}_{a_{1}}(x)  &   &    &    \\[0.5em]
  &    \mathcal{C}_{a_{2}(x)} &    &      &       \\[0.5em] 
&    &    \ddots    &  &      \\[0.5em]
&       &     &  &        \mathcal{C}_{a_{k}(x)}      
\end{pmatrix} 
\end{equation}

\begin{theorem} \label{RCFEquivSNF} Let $ M $ be an $ R[x] $-module that is free of finite rank $ m $ over $ R$ and  $ \beta $ be an $ R$-basis for $ M $.  Then $M $ admits an invariant factor decomposition if and only if the matrix   $   B(x) =  x I_{m}  -   [    \varphi ] _{     \beta }     $ has Smith normal form over $ R[x] $. 
\end{theorem} 

\begin{proof} $( \! \implies \! ) $ By the discussion above, we can find an $ R$-basis $ \alpha $ for $ M $ such that   $[\varphi]_{\alpha} $ is as in (\ref{xalpha}).  Let $ A(x) = xI_{m} - [\varphi]_{\alpha} $. It suffices to show that $ A(x) $ has SNF over $ R[x] $.  Indeed, if $ S \in \GL_{n}(R)  $ denotes the change of  basis   matrix from $ \beta $ to $ \alpha $  and if $ P(x), Q(x) \in \GL_{m}(R[x]) $ are such  that  $ P(x)   \cdot A(x)  \cdot Q(x) $ is in SNF,  then $ P(x) S  \cdot B(x) \cdot S^{-1} Q(x) $ returns  the same SNF.  Now  
 $$  A(x)  = \begin{pmatrix}   xI_{d_{1}}  -  \mathcal{C}_{a_{1}}(x)  &   &    &    \\[0.2em]
  &   x I_{ d_{2}}  -    \mathcal{C}_{a_{2}(x)} &    &      &       \\[0.2em]
&    &    \ddots    &  &      \\[0.2em]
&                       &                                      &  &       xI_{d_{k}}  -   \mathcal{C}_{a_{k}(x)}      
\end{pmatrix}    . $$   
It is easily seen by elementary row and  column operations that $ \mathrm{SNF}\big (xI_{d_{i}} -  \mathcal{C}_{a_{i}}(x)\big ) = \mathrm{diag} (1,\ldots,1,a_{i}(x)) $. Thus the matrix above is congruent over $ R[x] $ to $$ \diag( 1,  \ldots,   1   ,  a_{1}(x) , 1 , \ldots, 1 ,  a_{2}(x) , \ldots, 1, \ldots,  1 ,   a_ { k }     ( x  )   ) . $$   
Applying a few elementary row and column operations of switching, one obtains the desired SNF. \\

\noindent $( \! \impliedby \! )$ Say $ \beta = (b_{1}, \ldots, b_{m}) $.   Let $    \mathcal{ M } : = R[x] ^{  m } $ and  $ e_{1}, \ldots, e_{m}  $ be the standard $ R[x] $-basis for $ \mathcal{M}   $. Let $ \psi  : \mathcal{  M   }  \to  M $ be the surjective  homomorphism that sends $ e_{i} $ to $ b_{i}  $. It induces an isomorphism  $     \mathcal {  M   }    /  \ker \psi        \cong  M $. Suppose that $   [   \varphi             ] _{    \beta  } =      (  a_{i,j} )  $. For each  $ j \in \left \{ 1, \ldots, m \right \} $, set   $$    u_{j}  =  x e_{j} -  \sum_{ i = 1 }^{m } a_{i, j } e_{i    }     \in \mathcal{M} . $$ 
Then $ u_{j} \in \ker  \psi      $. Since $ x e_{j} - u_{j} $ lies in the $ R $-span of $ e_{1}, \ldots, e_{m} $, we see that 
\begin{align*}
 \mathcal { M  } &   =  R[x] e_{1}  +   \cdots + R [x]e_{m} \\
 &  =    R [x]  u_{1}  + \cdots  + R [x]u_{m} +   R e_{1} + \cdots +  R  e_{m}. 
 \end{align*}  
In other words,  every  element of $ \mathcal{M}  $ can be written as a sum of an element in the $ R[x]$-submodule $ U $ generated by $  u_{1} , u_{2} , \ldots , u_{m} $,  and an element of the $ R $-submodule $ V $ generated by $ e_{1} , e_{2} , \ldots , e_{m}  $ inside $ \mathcal{M} $.  We wish to show that $ \ker \psi = U $.  Pick $   a   \in   
\ker \psi  $. Then, $ a  = u + v $ for some $ u \in U $, $ v \in V $. Then $  \psi(v) = \psi(u+v) = \psi(a) = 0  $. 
But $ \psi(v) $ is $ 0 $ if and only if $ v = 0 $, since  $ \beta $ is a basis for $ M $. So $ a = u  \in  U  $. 

Consider now the matrix $ B(x) $. Its  $ j$-th column is the coordinate vector of $ u_{j} $ computed w.r.t.\  the standard basis $ ( e_{i} )_{i} $ of  $ \mathcal {  M  } $. Therefore right multiplication of $ B(x) $ by elements of $ \GL_{m}(R[x]) $ amounts to changing the set of generators of $ \ker \psi      $    and left multiplication amounts to changing the  standard  basis of $   \mathcal{M} $. Thus saying that $ B(x) $ has SNF over $ R[x] $ is equivalent to saying that there is an $R[x]$-basis of $    \mathcal  {  M   } $, say $ f_{1}, f_{2}, \ldots, f_{m} $    and polynomials $ a_{i} \in R[x]  $ satisfying $ a_{1} | a_{2 }  | \ldots | a_{m}  $ such that $  a_{1} f_{1}, \ldots , a_{m}  f_{m} $ forms a set of generators   of $ \ker \psi $.  This immediately implies that $ \ker \psi $ is a free $R[x] $-module and        
\begin{align*}   M    &  \cong     \mathcal{M }      /  \ker \psi    \\     & \cong  \bigg ( \bigoplus_{j=1}^{m} R[x] f_{i}    \bigg ) \bigg /   \bigg ( \bigoplus_{i=1}^{m} R[x] a_{i} f_{i}  \bigg  )   \\ 
& \cong \bigoplus_{i=1}^{m}  R[x]/ (a_{i})                                                                           
 \end{align*} 
All of the $ a_{i} $ can be taken to be monic, since the product of $ a_{i} $ equals $ \det(A(x)) $ which is monic.  If some $ a_{i} $ is $ 1 $, we can discard  the corresponding summand in the direct sum  decomposition above since  $  R[x] /  (a_{i}) = 0    $ in that case.                                                                                    \end{proof}

\begin{remark} The proof of the the backward direction  is  inspired by exercises 22-25 of \cite[\S 2]{D&F}.   
\end{remark}       
 
\begin{definition} We say that $ A \in \mathrm{Mat}_{m \times m }(R) $ has \emph{canonical form over $ R $}     if its rational canonical form (when $ A $ is seen over $ F$) lies in $ \mathrm{Mat}_{m \times m }(R) $ and is a $ \GL_{m}(R) $-conjugate of $ A $. In the case $ R = \ZZ $, such a form is also  referred to as an \emph{integral canonical form}.    
We similarly define this notion for endomorphisms of a  free $ R  $-module of finite rank.   
 \end{definition}                        
It is easy to see that an $ R[x]$-module $ M $ that is free of finite rank over $ R $ admits an invariant factor decomposition if and only if the map $ \varphi $ induced by $ x $ 
has canonical form over $ R $.    By replacing $ x    $ with $ -x $ in the theorem above, we obtain  the following.                                                                                                                        

\begin{corollary}              \label{RCFEquivMiller}                      For  an  $ r $-differential poset, Conjecture \ref{mainconj} is true for some integer  $ n $ if and only if 
$ DU_{n} : \ZZ P_{n} \to \ZZ  P_{n}    $ 
has  integral  canonical form.     
\end{corollary}

\section{Lifting canonical forms}  

For all of this section, $  R $ denotes  a  principal ideal domain. For any homomorphism $ \psi $ between two free $R$-modules of finite rank,  the rank of $  \psi $ equals the number of non-zero entries in a Smith normal form for $ \psi $.  By analogy, we refer to the number of unit entries as the  \emph{unital  rank} of $ \psi $. It describes the  size of any maximal subset  of  $ \mathrm{im}(\psi) $ that can be extended to a basis of the target of  $ \psi $. The nullity of $ \psi $ is denoted by $ \mathrm{null}(\psi) $.       
If $ A $ is a free $R$-module of finite rank, an element $ x \in A $ is said to be \emph{primitive} if $ x \notin IA $ for any proper ideal $ I \triangleleft A $.  Then primitive elements are precisely those that can be extended to an $ R$-basis for   $A $.   

\subsection{Tweaking  decompositions}   Fix now  an $R[x]$-module $ M$  that admits an invariant factor decomposition and let $ a_{1}, \ldots, a_{k} $ denote  its invariant factor sequence.    
\begin{lemma} \label{tweak}    Suppose $ M =  M_{1} \oplus \cdots \oplus M_{k} $ is an invariant factor decomposition and $ k \geq 2 $.   Let $ j , \ell  \in \left \{1, \ldots , k \right \} $ be two distinct  indices  and $ c \in R[x] $ be arbitrary.  Set $$   v_{\ell} '  :    =  \begin{cases} v_{\ell} + c v_{j} & \text{ if } j < \ell \\[0.5em]    v_{\ell}  + \dfrac{ c  a_{j}}   {  a_{\ell}    }      v_{j}  & \text{ if } j > \ell 
\end{cases}  
$$  
and let $ M_{\ell} ' $ be $ R[x] $-submodule of  $ M $  generated by $ v_{\ell} ' $. Then $ M_{1} \oplus \cdots \oplus M_{\ell}' \oplus \cdots \oplus M_{k}  $ is also an  invariant factor decomposition for $ M $.    
\end{lemma}

\begin{proof}  Since $ v_{\ell}' $ is an $R[x]$-linear combination of $ v_{\ell} $ and $ v_{j} $ with the coefficient of $ v_{\ell} $ being $ 1 $ and $ S =  \left \{  v_{1},  \ldots,  v_{k}   \right \} $ is a generating set for $ M $, the set $  \left \{ v_{\ell} ' \right \} \cup ( S \setminus \left \{ v_{\ell} \right \} ) $  generates $ M $  as  well.     Say $ p_{1} ,  \ldots,  p_{k} \in  R[x]  $ are such that
\begin{equation}  \label{summandvi} \displaystyle p_{\ell} v_{\ell} ' +  \sum_{ i \neq \ell} p_{i} v_{i}  =  0  
\end{equation}      We divide into two cases. \\[1em]  
\noindent \textit{Case 1: $j < \ell $.}\\[0.3em]  
We can rewrite (\ref{summandvi}) as $   (p_{j} + c p_{\ell}  )   v_{j  }    +   \sum _{ i \neq j }  p_{i}   v_{i}   =  0   $. As $ M_{1} , \ldots, M_{k} $ form a direct sum,  $ p_{i} v_{i} = 0 $ for all $ i \neq j $ and   $ (p_{j} + cp_{\ell}) v_{j}  = 0  $.  Since $ a_{i} $ generates  $  \Ann  (M _{i })    $ for all $ i $, we see that   $$  a_{j}  \mid ( p_{j}  + c p_{
\ell}) \quad  \text{ and }  \quad    a_{\ell}  |  p_{\ell}  . $$ 
Since $ a_{j} | a_{\ell} $ and $ a_{\ell} | p_{\ell} $, $ p_{\ell} $ annihilates both $ v_{j} $ and $ v_{\ell} $. Hence, it annihilates $ v_{\ell} ' = v_{\ell} + c v_{j} $. Similarly since $ a_{j} $ divides  both $ p_{j} + c p_{\ell} $ and $ p_{\ell} $, it divides $ p_{j} $ and therefore $ p_{j} $ annihilates $ v_{j} $.  Thus all summands in (\ref{summandvi})  vanish, and so the modules $ M_{i} $ for $ i \neq \ell $ 
form a direct sum with $ M_{\ell} '  $. Clearly $ a_{\ell } $ annihilates $ M_{\ell} ' $. Moreover any  
$ f \in \mathrm{Ann}(M_{\ell}  '    )  $ satisfies $  f v_{\ell} +  f c v_{j} = 0 $ which by the direct sum property  implies that   $ f v_{\ell} = fcv_{j} = 0 $. So $ f \in ( a_{\ell} ) $ and thus  $ (a_{\ell} ) =  \mathrm{Ann}(M_{\ell}') $.    
\\[1em] 
\noindent  \textit{Case 2:   $ j > \ell $.}\\[0.3em]    
 In this case, (\ref{summandvi}) can be rewritten as $ (p_{j} + c p_{\ell}  a_{j} /  a_{\ell}  ) v_{j} + \sum_{i \neq j} p_{i} v_{i} = 0 $.  As before, this implies that $ p_{i} v_{i} = 0 $ for all $i \neq j $ and $ (p_{j} + c  p_{\ell} a_{j} / a_{\ell}  )  v_{\ell} = 0 $.  Thus $$  a_{j} |  (  p_{j} + c p_{\ell} a_{j} / a_{\ell} ) \quad \text{ and } \quad a_{\ell} | p_{\ell} . $$ 
So $     p_{\ell} v_{\ell} '  = p_{\ell}    v_{\ell} +  c  \left (  p_ { \ell}  /   {  a_{\ell}  }  \right )    a_{j} v_{j}         =  0 $ and $    p_{j} v_{j}  =     (  p_{j} + c p_{\ell}  a_{j}/ a_{\ell}  )   v_{j}  -    c ( p_{\ell} / a_{\ell} )  a_{j} v_{j} = 0  $. This  establishes the direct sum property.  One similarly verifies  that $ a_{\ell} $ generates $ \mathrm{Ann}(M_{\ell} ' ) $.          
\end{proof}


\subsection{Ascensions}      \label{ascensec}    Let $ M $ be as above 
and $ N $ be an arbitrary $R[x]$-module that is free of finite rank over $ R  $. Suppose that    
$$ U : M \to N \quad \quad  \quad \quad D : N \to M  $$ 
are $ R[x]$-module homomorphisms such that $ DU : M \to M  $, $  UD : N \to N $ coincide with multiplication by $ x $, $ D $ is surjective, $ DU $ is injective and $ \mathrm{coker}(U) $ is free over $ R $.   

\begin{definition} An \emph{ascension} for an invariant factor decomposition $ M_{1} \oplus \cdots \oplus M_{k}  $ of $ M $ is a sequence $ W_{1} ,   \ldots,  W_{k} $ of monogenic $ R[x]$-submodules of $ N $ such that $ D(W_{i})  =  M_{i}  $.  We say that  the ascension is \emph{split} if   $ W_{1} , \ldots, W_{k} $ form a direct sum.  An \emph{ascension for $ M $}     is an ascension for some invariant factor decomposition of $ M $.    
\end{definition} 

Let $ W_{1}, \ldots, W_{k} $ be an ascension for $ M $. If $ w_{i} $ is a generator for $ W_{i} $, the image $ v_{i} = D(w_{i}) $ is a generator for $ 
D(W_{i}    )     $. So an ascension for $ M = M_{1} \oplus \cdots \oplus M_{k} $ can equivalently be described as a choice of generators $ v_{i} \in M_{i} $ for each $ i $ and a choice of $ w_{i} \in D^{-1} (v_{i}) $, though different choices may yield the same ascension. 

\begin{lemma}  \label{annascen}     The annihilator of the $i$-th member in an ascension for $ M $   is    either $ (a_{i})  $ or $ (x a_{i}) $. 
\end{lemma} 

\begin{proof}   Let $ W_{1}, \ldots , W_{k} $ be an ascension $ M _{i} : = D(W_{i}) $,  $ w_{i} \in W_{i} $ be generators and $ v_{i} = D(w_{i}) $.     First note that  $ x W_{i} = U(D (W_{i})) = U(M_{i}) $. Since $ DU $ is injective, so is $ U$. Hence $ U(M_{i})  $ is isomorphic to $ M_{i} $ which means that $  \mathrm{Ann}(U(M_{i})) = (a_{i}) $. Therefore $ xa_{i} $ annihilates $ W_{i} $.  Now either $ x a_{i} $ generates $ \Ann(W_{i}) $ or it does not.     If latter, pick any  $ p \in R[x] \setminus (x a_{i} ) $ such that $ p W_{i} = 0 $.  Then $ 0  =  D(p W_{i}) = p M_{i}  $  forces  $ p \in (a_{i} ) $. So $ p = c a_{i} $ for some $ c \in R[x] $. Since $ x a_{i} $ annihilates $ W_{i} $,  we may assume wlog that $ c $ is a non-zero element of $ R $.  But since $  N $ is free over $ R $, $ ca_{i} W_{i} = 0 $ implies $ a_{i}W_{i} = 0 $. So $  \mathrm{Ann}(W_{i} ) = (a_{i}) $ in this  case.  
\end{proof}

\begin{proposition}    \label{ascensionisdirect}   If $ a_{1}(0) \in R $ is non-invertible, any ascension for $ M $  is  split  and the annihilator of the $ i $-th  ascended  member equals $ ( x a_{i} ) $.  
\end{proposition}      
\begin{proof} 
 Let $ W_{1}, \ldots , W_{k} $ be an ascension $ M _{i} : = D(W_{i}) $,  $ w_{i} \in W_{i} $ be generators and $ v_{i} = D(w_{i}) $.    
Suppose that $ p_{1}, \ldots, p_{k} \in R[x] $ are such that 
\begin{equation} \label{dirsumtest}  p_{1} w_{1} +  \cdots + p_{k} w_{k} = 0  \end{equation}  
Applying $ D $ gives $  p_{1} v_{1} + \cdots + p_{k} v_{k} = 0 $. Therefore $ p_{i}  \in ( a_{i} ) $ for all $ i = 1, \ldots , k $. Let $ c_{i} \in R[x] $ be such that $ p_{i} = c_{i} a_{i} $. Since $ xa_{i} $ annihilates $ w_{i}  $ by Lemma \ref{annascen},  we may assume wlog that $ c_{i}   \in  R $.  Let $ I = (c_{1}, \ldots, c_{k}) \subset R $ be the ideal generated by these.  If $ I = (0) $, all of the $ p_{i} $ vanish, so assume otherwise. By replacing $ p_{i} $ with $  \gamma  ^{-1} p_{i} $ for a generator $ \gamma \in I $, we may also assume wlog  that    $ I = R $. Now (\ref{dirsumtest}) can be rewritten as
\begin{equation}  \label{dirsumtestconseq} \sum_{i=1}^{k}  c_{i} \big ( a_{i}(x) - a_{i}(0)  \big  )w_{i}  =   -  \sum_{i=1}^{k} c_{i} a_{i}(0) w_{i} 
\end{equation}   
Let $  \omega    : =  \bigcup _{ i =1 }^{k} \left \{ x^{j} w_{i} \, | \, 0 \leq j \leq d_{i}  -   1    \right \}  \subset N$ and $ \omega_{ 1 }    :   = x \omega  $.  Since $ DU $ acts as $ x $ and $ U $ is injective, we see that $ \omega  _   {  1  }   $ is an $ R $-basis for the image of $ U $ (see the diagram below).
\bigskip 
\begin{center}
\begin{tikzcd}[ 
        ,every arrow/.append style={maps to}
        ,every label/.append style={font = \scriptsize}
] 
  
{w}_{i}  \arrow{r}\arrow{dr}[swap]{D}   &   [0.6cm]     x {w}_{i} \arrow{r}\arrow{dr}[swap]{D} &[0.6cm]    x^{2}{w}_{i}\arrow{r}\arrow[maps to]{dr}[swap]{D } & [0.6cm]x^{3}  {w}_{i}\arrow{r}\arrow{dr}[swap]{D}&  [0.6cm] \cdots\arrow{r}\arrow{dr}[swap]{D }    &  [0.6cm]  x^{d_{i}}  {w}_{i}   \\ [0.6cm]
&     \, \,  v_{i}     \,    \,       \arrow{u}[]{U}\arrow{r}[]{}&      \, \,      x v_{i}  \,  \,   \arrow{r}\arrow{u}{U}&     \, x ^{2}  v_{i}\,   \arrow{u}{U} \arrow{r} & \cdots                 \arrow{r}     &                 x ^ { d_{i} - 1}v_{i}       \arrow{u}{U}   
\end{tikzcd}
\end{center}        
\bigskip  
Since $  \mathrm{coker}(U) $ is free over $ R $, the image $ \mathrm{im}(U) = R \omega_{1}    $ is a direct summand of $ N $, i.e.,    $ \omega_{1} $ can be extended to an $R$-basis for $ N $. Now the LHS of (\ref{dirsumtestconseq}) is an $ R$-linear combination over $  \omega_{1} $ and the monicity of the $ a_{i} $ implies that the ideal generated by the coefficients of this linear combination equals $ R $.   
Thus the LHS of  (\ref{dirsumtestconseq}) is  a primitive 
element of $ N $. 
The RHS however is an element of $ a_{i}(0) N $ since $ a_{i}(0) \in ( a_{1}(0)) $ for all $ i $. As  $ a_{1}(0)   \notin R^{\times} $, this contradicts the primitivity of the LHS.  Thus $ I $ must vanish and so must    $ p_{i} $.  Therefore $ W_{1}  + \cdots  +  W_{k }    $ is a direct sum. The same argument reveals that $ p \in R[x] $ annihilates $ w_{i}  $ if and only if $ p   \in (x a_{i} )$ which establishes the second claim. 
\end{proof}

\begin{corollary} \label{nullityofD}     If $ a_{1}(0) $ is non-invertible, $ \mathrm{null}(D) \geq k   $.    
\end{corollary}  
\begin{proof} This follows by noting that  $  a_{i}w_{i} \in W_{i} \cap \ker(D) $ and  $  \left \{ a_{1}w_{1} , \ldots, a_{k} w_{k} \right \} $ is an $R$-linearly independent subset of $ \ker (D) $ of  cardinality $ k  $. 
\end{proof}  
Given an ascension $ W_{1}, \ldots, W_{k} $, one can construct new ones by replacing the generators $ w_{i}  \in W_{i}  $  using  one  of  the following operations: 
\begin{enumerate}
\item[(A)] replace $ w_{\ell}  $ with $ w_{\ell} +  t$  for any index $ \ell $ and $ t \in \ker (D) $,   
\item[(B)]  replace $ w_{\ell} $ with $ w_{\ell} + c w_{j} $ for any $ j < \ell $ and  $ c \in R $
\item [(C)] replace $ w_{\ell} $ with $ w_{\ell} + c (a_{j}/a_{\ell} )w_{j}   $ for any  $ j  > \ell $ and  $ c \in R $.    
\end{enumerate} 
Operation (A) is justified since $ w_{\ell  }    + t  \in D^{-1}(v_{\ell} )$. Operations (B), (C) are justified by Lemma \ref{tweak} and amount to changing the underlying invariant factor decomposition for $ M $. 
 \begin{theorem}   \label{mainabsteo1}    If $ a_{1}(0) $ is non-invertible  and  $ \mathrm{null}(D) \geq   k + 1 $,  
$ N $ has an invariant factor decomposition.  
\end{theorem}  
\begin{proof}  Fix an ascension $ W_{1}, \ldots, W_{k} $ for $ M $ and let $ M_{i} $, $ w_{i} $, $ v_{i} $ be as above. 
Let $ \kappa_{0} :  = \left \{ a_{1}w_{1} , \ldots, a_{k} w_{k} \right \}   \subset \ker(D)   $     
and $ p :  = \mathrm{null}(D) $.  
Since $ R $ is a PID, the inclusion $ R \kappa_{0}  \hookrightarrow \ker D $ admits a Hermite normal form (cf.\  \cite[\S 3.2]{Zassenhaus}).  More precisely, there is an $ R $-basis $  t_{1}, \ldots, t_{p} $ of $ \ker ( D  ) $  such  that  for all $ j $,    
\begin{align}  \label{hermiteeq}  a_{j}  w_{j}     & =  \sum_{i = 1 } ^ {  j  }   c_{i,j} t_{i}    
\end{align}
where   $     c_{i,j} \in R $. In other words, the matrix of the inclusion $ R \kappa_{0} \hookrightarrow \ker (D) $ is $$ \quad \quad \quad  \quad  \quad  E =   \left (  \begin{matrix} c_{1,1} & c_{1,2} & \cdots & c_{1,k}   \\
& c_{2,2} & \cdots & c_{2,k} \\
 &  &  \ddots & \vdots  \\
 &  & & c_{k,k} \\
 & 
\end{matrix}  \right   )  \in  \mathrm{Mat}_{p \times k } (  R    )     .  $$
Let $ E_{i} $ denote the $ i $-th column of $ E $. Let us record how $ E $ changes when we modify the generators $ w_{i}   $ with respect to operations (A), (B), (C): 
\begin{enumerate} 
\item [(A)] replace $E_{\ell}  $ with $ E_{\ell} + a_{\ell} (0) F $  for any $ \ell $ and $ F  \in \mathrm{Mat}_{p \times 1}(R) $,    
\item[(B)]   replace $ E_{\ell} $ with $ E_{\ell} + c \big((a_{\ell}/a_{j})(0)\big) E_{j}  $  for any $ j   <  \ell $ and $ c \in  R $,    
\item [(C)] replace $ E_{\ell} $ with $ E_{\ell} + c E_{j} $ for any $  j    >  \ell  $ and $ c \in R $.        
\end{enumerate} 
We say that $ E $ is \emph{connected} to $ E ' \in \mathrm{Mat}_{p \times k}(R) $ if there exists  $ P \in \GL_{p}(R) $ such that $ E' $ can be obtained from $ PE $ by applying operations (A), (B), (C) a finite number of times.  \\[1em]   
\noindent \textit{Claim:  $ E $ is connected to the matrix that has $ 1_{R} $ on its diagonal and zero elsewhere.} \\
Begin by placing $ a_{1}(0) $ in $ (1,2) $-position of $ E $ via operation (A) and call this matrix $ E '$ for the moment.  We claim that ideal $ J $  generated by $ a_{1} (0) $ and $ c_{1,1} $ equals $ R $. Suppose not.   
Then  (\ref{hermiteeq})  implies that $$ \left ( a_{1}(x) - a_{1}(0) \right ) w_{1} = c_{1,1} t_{1}  -  a_{1}(0) w_{1}  \in J N . $$ 
This is however a contradiction, since the LHS is a primitive element of $ N $ (see the proof of Proposition \ref{ascensionisdirect}). So the first column of $ E $ describes a primitive vector in $ R^{2} \subset  R^{p} $ where ($ R^{2} $ is embedded in the first two components).   There is thus an element $ P \in  \GL_{p}(R) $ which is a block diagonal sum of a $ 2 \times 2 $ invertible matrix with the $ (k-2) \times (k-2)$ identity matrix such that the first column of $ PE' $ is $ (1,0, \ldots, 0)^{t} $. Operation (C) then allows us to replace $ PE' $ with a matrix $ E'' $ whose first row is $ (1, 0, \ldots, 0) $. Relabeling everything, we can assume that $$ E =      \left (  \begin{matrix} 1  &   & 
&  \\
& c_{2,2} & \cdots & c_{2,k} \\
 &  &  \ddots & \vdots  \\
 &  & & c_{k,k} \\
 & 
\end{matrix}  \right   )  . $$
Since $ a_{1}(0) | \cdots | a_{k}(0) $, we can continue inductively:   place $ a_{i}(0) $ in the $ (i,i+1)$-entry, use primitivity of elements to establish $ (a_{i}(0), c_{i,i}) = R $, transform by multiplying on the left by an element of $ \GL_{2}(R) $ to make $ (i,i)$-entry $ 1 $ and then use this entry to make all entries to its right $ 0 $. Note that $ p > k $ is necessary to execute this procedure in the $ k $-th  column.    The claim  is   established.   \\[1em] 
Thus we may assume that the ascension $ W_{1} ,  \ldots , W_{k } $ is such that $ R \kappa_{0} = \bigoplus_{i=1}^{k} R a_{i}w_{i}  $ is a direct summand of $ \ker(D) $. Let 
$ \omega  : =  \bigcup_{i=1}^{k} \left \{  x^{j} w_{i} \, | \, 0 \leq j \leq d_{i} -1 \right \} $. Since $ D $ maps $\omega $ bijectively onto the $ R$-basis $ \left \{ x^{j} v_{i} \, | \, 0 \leq j \leq d_{i}-1     \right \} $ of $ M $,  we have $$ N =  R \omega  \oplus \ker (D)   . $$ 
Choose $ \kappa $ a basis of $ \ker(D) $ that extends $ \kappa_{0} $, so that $ \omega  \cup \kappa  $ is an $ R $-basis for $ N $. Then the set $ \nu  $ obtained by replacing each $ a_{i} w_{i} \in \omega \cup \kappa $ with $  a_{i} w_{i} + ( x^{d_{i}}     - a_{i})w_{i}  = x^{d_{i}} w_{i} $ is also an $ R $-basis since $ (x^{d_{i}} - a_{i})   w_{i} \in  R \omega $.  We can write $ \nu  $ as the union of  $ \bigcup_{i=1}^{k} \left \{ x^{j} w_{i} \,  |  \, 0 \leq j \leq d_{i}  \right \} $ with $ \kappa \setminus \kappa_{0} \subset \ker(D)  $. Then    if    $ s_{1}, \dots, s_{p -k } $ denote the $  p - k $ elements of $ \kappa \setminus \kappa_{0} $, $$  N = R \nu  =  R s_{1} \oplus \cdots \oplus R s_{p - k } \oplus  W_{1} \oplus \cdots \oplus W_{k}  . 
$$ 
is the desired  invariant factor decomposition for $ N $.   
\end{proof}      
\begin{remark} Note that the proof above goes through without ever invoking operation (B). However we chose to include it for the sake of completeness, as this operation seemed natural in the formulation of Lemma \ref{tweak}.     
\end{remark}

\subsection{Invertible constants} 

It will be necessary to generalize  Theorem \ref{mainabsteo1} to some cases where $ a_{1}(0) \in R^{\times} $.  We record such a generalization  below.       

We continue to assume the notation and terminology introduced so far.  Observe that the injectivity of $ DU $ implies that  the intersection $ \ker(D) \cap  \mathrm{im}(U) $ is trivial. 
Let $ \iota $ denote the inclusion $ \ker(D) \oplus \mathrm{im}(U)  \hookrightarrow N $. Note that this is a full rank inclusion. %
Let 
$ k_{0}   $ denote the  number of non-invertible terms in the sequence $ (  a_{1}(0) , \ldots,  a_{k}(0)) $.   
\begin{lemma} \label{smithofD}   A Smith normal form for the inclusion  $ \iota $ is given  by   $  \mathrm{diag}(1,\ldots,1, a_{1}(0), \ldots,  a_{k}(0)) $. In particular, $ \mathrm{null}(D)  \geq  k_{0}    $. \end{lemma}    
\begin{proof}      Let $ L : =  \ker(D) \oplus \im(U  )$. Then $ L $ is the kernel of  the composition  $   N   \xrightarrow{D }  M  \xrightarrow{\, \pi  \, }   M / x M  $ 
where $ \pi $ denotes the canonical quotient map. So we have an isomorphism $$ N /  L  \simeq  M / x M  \simeq \bigoplus _{ i  = k_{0} } ^   {  k  }   R/(a_{i}(0)).  $$
of $ R $-modules and the  first claim easily follows from this.    
Since  $ \mathrm{im}(U) $ is a  direct   summand of $ N $, the unital rank of $ \iota $ is at least $ \mathrm{rank}_{R}(M) $. Since this  unital rank 
is  exactly   $ \mathrm{rank}_{R} (N )  -  k_{0}    $ in light of  the first claim,  the second  claim also  follows.    
\end{proof}  

\begin{theorem}  \label{mainabsteo2}  Suppose that there  is  a positive integer  $ \delta   \leq  k    $ such that $ a_{1} = \cdots = a_{\delta} = x-1 $, $ a_{i}(0) $ is non-invertible for $ i > \delta $     and $  \mathrm{null}   (D) \geq k - \delta + 1 $.    Then $ N $ admits an invariant  factor   decomposition.    
\end{theorem}    

\begin{proof} Fix an ascension $ W_{1}, \ldots, W_{k} $ for $ M $ and let $ M_{i}, w_{i}, v_{i} $ be as  usual.  By replacing $ w_{i} $ with $ x w_{i} = w_{i} + (x-1)w_{i}  $, we may assume that $ x w_{i} = w_{i}  $ for $ i = 1, \ldots, \delta $. Then $ \mathrm{Ann}(W_{i}) = (x-1) $ for $1 \leq  i \leq \delta $ and the argument of Proposition \ref{ascensionisdirect} still goes through  to show that any such  ascension   splits and that $ \mathrm{Ann}(W_{i}) = ( x a_{i} ) $ for $ i > \delta $.    
Similarly the argument of Theorem \ref{mainabsteo1} can be executed to show that there is a modification of  the $   k - \delta  =  k_{0}     $ generators $ w_{\delta+1} , \ldots, w_{k} $ by elements of $ \ker(D) $ so that  $ \kappa_{0}  =  \left \{ a_{i} w_{i} \, | \, \delta < i \leq k \right \} $ extends to an $R $-basis of $ \ker(D) $. 
Let $  \omega  = \bigcup_{i=1}^{k} \left \{  x^{j}   w_{i} \, | \, 0 \leq j \leq d_{i} - 1  \right \}  $ as before and let $ \kappa  $  be  a   basis of $ \ker (D ) $ that extends $ \kappa_{0} $. Then $ \omega \cup  \kappa $ is an $ R $-basis for $ N $ and therefore so is $$ \nu '    : = \omega   \cup ( \kappa \setminus \kappa_{0  } )  \cup \left \{ x^{d_{i}    }  w_{i} \, | \, \delta < i \leq k \right \} . $$   Let $ p  =  \mathrm{null} (D)   $ and  $ s_{1}, \ldots, s_{p-k_{0}} $ denote the $ p - k_{0} $  elements of $ \kappa \setminus \kappa_{0} $. Note that $ p - k_{0} = p - k + \delta \geq 1 $.  
 \\[1em] 
\noindent \textit{Case 1: $ p - k_{0} \geq \delta $.} \\[0.3em]  
 Construct the $ R $-basis $ \nu : = \left \{  s_{i} +  w_{i} \,  | \,  1 \leq i \leq \delta 
  \right \} \cup ( \nu' \setminus \left \{  s_{1}, \ldots, s_{\delta} \right \}   )   $ for $ N $. If $ W_{i}^{\circ} $ denotes the $ R[x] $-submodule generated by $ w_{i} + s_{i} $ for $ i = 1 , \ldots, \delta $, then $ \mathrm{Ann}(W_{i}^{\circ}) = x(x-1) $. Therefore  
 $$ N = R \nu =  \left (    R s_{\delta + 1 } \oplus \cdots \oplus R s_{p-k_{0}}  \right   ) \oplus \left (   W_{1}^{\circ} \oplus \ldots  \oplus  W_{\delta}^{\circ}  \right ) \oplus \left ( W_{\delta + 1 } \oplus \ldots      \oplus W_{k} \right )   $$
is an invariant factor decomposition for  $ N $. \\[1em] 
\noindent \textit{Case 2: $p-k_{0} < \delta $}\\[0.3em] 
Construct the $ R$-basis $ \nu : =  \left \{  s_{i} + w_{\delta - i + 1   }  \, | \, 1 \leq i \leq   p - k_{0}    \right \}     \cup ( \nu' \setminus \left \{  s_{1}, \ldots, s_{ p - k_{0}}   \right \}   )   $.  Let $ W_{\delta - i  + 1 }^{\circ} $ denote the $ R[x]$-submodule generated by $ s_{i} + w_{\delta - i + 1 }  $ for $ i = 1, \ldots, p -  k_{0} $. Then again, 
$$ N = R \nu = \left (  W_{1} \oplus \ldots \oplus W_{ \delta  + k _{0 } - p } \right )  \oplus  \left (  W_{ \delta  + k_{0}  - p 
   + 1 }^{\circ}   \oplus \ldots \oplus  W_{\delta } ^\circ \right )    \oplus \left (   W_{\delta+1} \oplus  \ldots \oplus  W_{k}  \right )    $$  
is an invariant factor decomposition.    
\end{proof}

\section{The main theorem}    
We now apply the result of the previous section to differential posets. Let $  P $ denote an $ r $-differential poset, $ P_{n} $ its $ n $-th rank (which is empty if the integer  $ n $ is negative),  $ p_{n} $ the cardinality of $ P_{n} $ and $ \Delta p_{n} $ the difference  $  p_{n} - p_{n-1} $. Define $  U_{n} $, $ D_{n} $ etc., as in the  introduction  for  $ R = \ZZ $. Let $ \delta_{i,j} $ denote the Kronecker delta.     

\begin{lemma}       \label{freecokernel}                For  an  $ r $-differential poset, the map $ D_{n+1} $ is surjective    if and only if $ U_{n}   $ has free cokernel.                          
\end{lemma}

\begin{proof} That $ U_{n} $ has free cokernel is equivalent to saying that the Smith normal form of $ U_{n} $ consists of only ones and  zeros.  The same applies to $ D_{n+1} $. But since $ D_{n+1} $ is of full rank, freeness of cokernel is equivalent to  surjectivity. Now note that the matrices of  $ D_{n+1} $, $ U_{n} $ in the standard basis are transposes  of each other and thus the same is true for their  Smith normal forms. 
\end{proof}

\begin{proposition} \label{conjequiv}  Suppose that $ D_{n+1} $ is surjective and  $ DU_{n}  $ has integral canonical form for some $ n $.    
Then $\Delta p_{n+1}  \geq \Delta p_{n-i-\delta_{r,1} }  $
for all $ i  \geq  0  $.    
\end{proposition}

\begin{proof} We proceed via  induction.  The case $ n = 0 $ is trivial. Say the claim is true for some $ n \geq 0 $.     
Set $$ M :   =  \ZZ P _{n }, \quad   \quad  \quad    N  :    =  \ZZ    P_{n+1} $$ and view these are $  \ZZ  [   x] $-modules via $ DU_{n} $, $ UD_{n+1} $ respectively.  Let  $ D : N \to N $, $ U : M \to N $ be the up and down maps respectively. Note that both of these respect the $\ZZ[x]$-module structures. Moreover, our induction hypothesis implies that $ M $ has invariant factor decomposition. Thus we are in the setup of \S \ref{ascensec}.  By the discussion in \S \ref{recoll}, the integer $ k $ for the invariant factor decomposition for $ M $ equals  $ \max \left \{ \Delta p_{0}, \ldots, \Delta p_{n}  \right \}  $. 
 If $ r >  1  $, Corollary \ref{nullityofD} and the expression (\ref{exp}) for $ a_{1}(x) $ implies that $ \Delta p_{n+1} = \mathrm{null}(D)  \geq k $.  If $ r = 1 $, the integer $ k_{0} = \max \left \{ \Delta p_{0}, \ldots, \Delta p_{n-1} \right \} $ computes the number of invariant factors of $ M $ whose constant term is non-invertible in $ \ZZ  $. In this case,  Lemma \ref{smithofD} implies the claim. 
 \end{proof}

 \begin{theorem}  \label{mainteo} Suppose that  for the $ r $-differential poset  $  P  $, 
\begin{enumerate}[itemsep=0.1em]   
\item[$\bullet$]
all the down maps  are surjective,    
\item[$\bullet$]
there exists a non-negative integer $ m  $ such that   $   \Delta   p_{n}     >  \Delta p_{n-1 - \delta_{r,1}}    
 $ for every $ n > m                       $, 
\item[$\bullet$] $ DU_{0}, \ldots, DU_{m}  $ have    integral canonical forms.  \end{enumerate}              
Then $ DU_{n}  $ has integral canonical   form  for  every  $   n $.                       
\end{theorem} 
\begin{proof}  We  proceed via strong induction 
on $ n $. Base case verification is included in the third bullet. Say all of $  DU_{0}, \ldots, DU_{n} $  have integral canonical forms for some $ n \geq  m  $.   Denote as in Proposition \ref{conjequiv}  $ M = \ZZ P_{n} $, $ N = \ZZ P_{n+1}  $ and write  $ D $, $ U $ for up and down maps between $ M $, $ N $. Then $ M $ has invariant factor decomposition with $ k = \max \left \{ \Delta p_{0}, \ldots, \Delta p_{n} \right \} $ components. Proposition \ref{conjequiv} implies that $$ k = \begin{cases}   \Delta p_{n} & \text{ if }  r > 2 \\ \max \left \{ \Delta p_{n} , \Delta p_{n-1}  \right \}  & \text{ if } r = 1 . \end{cases} $$
  If either $ r = 2 $ or if $ r =1 $ and $ k = \Delta p_{n-1}  $, we see from (\ref{exp}) that $ a_{1}(0) $ is non-invertible in $ \ZZ $. So Theorem \ref{mainabsteo1}  implies   that $ N $ has an invariant factor  decomposition. If on the other hand 
  $ r = 1  $ and  $ k = \Delta p_{n}  > \Delta p_{n-1} $,  the integer $ \delta =  \Delta p_{n} - \Delta p_{n-1} \geq 1 $ is such that $ a_{1} = \ldots = a_{\delta} = x- 1 $ and $ a_{i}(0) $ is non-invertible for $ i > \delta $. In this case, Theorem \ref{mainabsteo2} implies that  $  N $ has an invariant factor decomposition. 
In either case, we see that  $  UD_{n+1} $ has integral canonical form. Equivalently, $  xI_{p_{n+1}}  - [UD_{n+1}] $ has Smith normal form over $ \ZZ[x] $ (see    Theorem \ref{RCFEquivSNF}).  Since $$ DU_{n+1} = UD_{n+1}  +  r   \cdot     \mathbbold{1} ,  $$   $ xI_{p_{n+1}} - [DU_{n+1}] $ also has Smith normal form over $ \ZZ[x] $ which is given by replacing $ x $ with $ x - r $ in the Smith normal form of $ UD_{n+1} $. Then again, this implies that $ DU_{n+1} $ has integral canonical form.  This completes the induction step.
\end{proof}

\begin{remark}   Combining Proposition \ref{conjequiv} with  \cite[Proposition 2.5]{Miller}, we see that Conjecture 2.3 and 2.4 of \emph{loc.cit.}\ are equivalent whenever Conjecture \ref{mainconj} holds. Then  Theorem \ref{mainteo} may  be seen as a converse of sorts.           
It would be interesting to relax the second condition of Theorem \ref{mainteo} to the  bound given in Proposition \ref{conjequiv}. 
\end{remark}

\section{Applications}        

In this section, we  record some  applications.    

\begin{theorem} \label{Cartesian}    Let $ P $ and $ Q $ be differential posets of  rank sizes $ p_{n} $, $ q_{n} $ respectively. Suppose that $ \Delta q_{n} \geq \Delta q_{n-1}    $ for all $ n \geq 2 $, and that all the down maps of at least one of the posets are surjective. Then  Conjecture \ref{mainconj} holds for  $ P  \times  Q    $. 
\end{theorem}

\begin{proof}     By Corollary \ref{RCFEquivMiller}, it is enough to show that the $ DU_{n} $ maps of $ P \times Q $ have integral canonical form.                                                                             Notice that $ P \times Q $ is an  $ r$-differential poset for some $ r \geq 2 $. It was proved in \cite[Proposition 4.5]{Miller} that the up maps of a cartesian product have free cokernel if one of the posets in the product has this property. So Lemma \ref{freecokernel}  implies that  the down maps of $ P \times Q $  are surjective. Denote the rank sizes of  $ P \times Q $ by $ \rho_{n} $. Then 
$     \rho _{n}  =                         \sum_{  i = 0 } ^{  n}      q_{n-i}                   p _{ i }               $, so 
\begin{equation} \label{rankexp}  \Delta \rho _{n} - \Delta \rho _{n-1}   =  q_{0}(    \Delta p_{n}     ) +   \Delta  q_{1}  p_{n-1}  +     \sum_{i=0}^{n-2}         ( \Delta  q_{n-i}  - \Delta q_{n-i -1 }  )    p_{i}                           .                             
\end{equation} 
Suppose that  $ n \geq 2 $. If $ Q $ is $ 1  $-differential, then $ \Delta q_{2} - \Delta q_{1} = 1$. So the last summand of the sum in (\ref{rankexp})                           contributes a non-zero term. If $ Q $ is $  s $-differential for some $ s  > 1 $, then 
$ \Delta q_{1} =  s - 1 $ which means that $ \Delta q_{1} p_{n-1} \geq 1  $.  Since all the terms in the sum (\ref{rankexp})   are non-negative,       we see that  $$ \Delta \rho _{n}   > \Delta \rho _{n-1}   $$ for all $ n \geq 2 $. 
Additionally if  $ PQ $ is $ r $-differential for $ r \geq 3 $, then $ \Delta \rho_{1} - \Delta \rho_{0}  = r - 2 $ is at least $ 1 $.  Now  $ DU_{0} $ trivially has integral canonical form, so Theorem \ref{mainteo} with $ m = 0 $  gives the result.  
If $ r = 2 $, the matrix for $ DU_{1} $ is always $ \left ( \begin{smallmatrix} 3 & 1 \\ 1 & 3     \end{smallmatrix} \right )    $, and one can easily verify that it has an integral canonical form. 
So in this case, 
Theorem \ref{mainteo} applies with $ m = 1 $.    
\end{proof}

\begin{corollary}            The Conjecture \ref{mainconj} is true for $ \mathbf{Y}^{r} $ for every $ r \geq 1  $.
\end{corollary}  
\begin{proof} Suppose first that $ r = 1 $.   It was proved in \cite[\S   6.1]{Miller} that the up maps of  $     \mathbf{Y}  $ have free cokernel and hence  the down maps are surjective. We claim that $$ \Delta p_{n}  >  \Delta p_{ n - 2} $$ holds for  all $ n > 2  $.  
Since $  \Delta  p_{1} =    0  $ and $ \Delta p_{3} = 1 $, the condition  holds for $ n = 3 $.    So assume $ n \geq   4 $.  Notice that $ \Delta p_{n} = p_{n} - p_{n-1} $ counts the number of partitions of $  n $ with no part equal to $ 1 $. Let $ S_{n} $ be the set of all such partitions     of $n $. For each partition in $ S_{n-2} $, we can add a $ 2 $ to the largest part, and obtain a partition of $  n    $    in                                                      $ S_{n} $. This injects $ S_{n-2} $ in $ S_{n} $ and so $ |S_{n} |  \geq  | S_{n-2} |.  $   If $ n $ is even, the partition $ 2, 2, 2, \ldots, 2 $ with $ n/2 $ number of $ 2 $s cannot be obtained from the said injection of $ S_{n-2} $ into $S_{n} $. Similarly if  $ n $ is odd, the partition $ 3, 2, 2, \ldots, 2      $ with $ \lfloor n/2 \rfloor    -   1      $ number of $ 2 $s does not arise from a partition in  $  S_{n-2} $.  So  we have $$ | S_{n} |  > |  S_{n-2}  |   $$ for $ n  \geq  4                                                             $, and we obtain the desired inequality. 
Now one can easily verify that $ DU_{0} $, $ DU_{1} $ and $ DU_{2} $ for $ \mathbf{Y} $  all have integral canonical form. Invoking Theorem \ref{mainteo} (with $ m =  2$)  and Corollary    \ref{RCFEquivMiller}, we get the result in this  case.

For the case  $ r > 1 $, note that there is an injection $ S_{n-1} $ in $ S_{n} $ given by adding $   1  $ to the largest part in a partition from $ S_{n-1} $. This  implies that   $ |S_{n}| \geq |S_{n-1}| $ for all $ n \geq 2 $. The claim then  follows by Theorem \ref{Cartesian} with $ P =  \mathbf{Y}^{r-1} $ and  $ Q = \mathbf{Y} $. 
\end{proof}

\begin{corollary}   The Conjecture \ref{mainconj}     holds for $ Z(r) $ for all $ r \geq 1 $.   \end{corollary}

\begin{proof} The surjectivity of down maps for $ Z(r) $  was proved in \cite[\S 5]{Miller}. Recall that the rank sizes of $ Z(r) $ satisfy the recursion given by $ p_{0} = 1 $, $ p_{1} = r $  and  $ p_{n} = r p_{n-1} + p_{n-2} $ for $ n \geq 2 $. Therefore the rank differences satisfy $ \Delta p_{0} = 1 $,  $ \Delta p_{1} =  r  -  1  $ and $$  \Delta p_{n} =   (r-1) ( \Delta p_{n-1} )  +  r  p_{n-2} ,  \quad \quad  n \geq 2 . 
$$ 
So if $ r > 2 $,  $ \Delta p_{n} > \Delta p_{n-1} $ holds  for all $ n > 0 $ and Theorem \ref{mainteo} applies with $ m = 0 $ (the base case is trivial). If $ r = 2 $, $ \Delta p_{n} > \Delta p_{n-1} $ for $ n > 1 $ and Theorem \ref{mainteo} applies with $ m = 1 $ (base case for $ DU_{1}$ is the same as in Theorem \ref{Cartesian}). Finally if $ r = 1 $, $ \Delta p_{n} > \Delta p_{n-2} $ holds  $ n  >  3   $. It is easily verified that $ DU_{0} $, $ DU_{1} $, $ DU_{2} $ and $ DU_{3} $ all have integral canonical forms.  Therefore Theorem \ref{mainteo} applies with $ m = 3 $.    
\end{proof}

\begin{remark} This result   was also proved in \cite[\S 5]{Miller}.  
\end{remark}

\bibliographystyle{amsalpha} 
\bibliography{references}{}

\end{document}